\documentclass[12pt,reqno,a4wide]{amsart}

\usepackage{tikz} 

\allowdisplaybreaks

\newtheorem{theorem}{Theorem}

\begin{document}

	\bibliographystyle{plain}

	\title[The generating function of $p$-Bernoulli numbers: An elementary proof]{A closed formula for the generating function of $p$-Bernoulli numbers: An elementary proof}

	\author[Helmut Prodinger]{Helmut Prodinger}
	\address{Department of Mathematics, University of Stellenbosch 7602,
		Stellenbosch, South Africa}
	\email{hproding@sun.ac.za}

	\author[Sarah J. Selkirk]{Sarah J. Selkirk}
	\address{Department of Mathematics, University of Stellenbosch 7602,
		Stellenbosch, South Africa}
	\email{sjselkirk@sun.ac.za}
	\thanks{The financial assistance of the National Research Foundation (NRF) towards this research is hereby acknowledged. Opinions expressed and conclusions arrived at, are those of the authors and are not necessarily to be attributed to the NRF}

	
	\begin{abstract}
		For a two parameter family of Bernoulli numbers $B_{n, p}$ the exponential generating function is derived by elementary methods. 
	\end{abstract}
	

	\maketitle

	\section{Introduction}
The following recursion for a two paramenter family of Bernoulli numbers is given in \cite{Rahmani}, 
\begin{equation}
\label{eqn:rec}
B_{n+1, p} = pB_{n, p}-\frac{(p+1)^{2}}{p+2}B_{n, p+1} \quad \text{for} \quad n, p \geq 0.
\end{equation}
In terms of exponential generating functions 
\begin{equation*}
f_{p}(t) := \sum_{n\geq 0}B_{n, p}\frac{t^{n}}{n!},
\end{equation*}
recursion (\ref{eqn:rec}) translates into 
\begin{equation*}
f_{p}^{\prime}(t) = pf_{p}(t)-\frac{(p+1)^{2}}{p+2}f_{p+1}(t).
\end{equation*}

The closed formula that follows is the main result of \cite{doi:10.2989/16073606.2017.1418762}.
\begin{theorem}
\label{thm}
For $p\geq 0$
\begin{equation}
\label{eq:main}
f_{p}(t) = \sum_{n=0}^{\infty}B_{n, p}\frac{t^{n}}{n!} = \frac{(p+1)(t-H_{p})e^{pt}}{(e^{t}-1)^{p+1}}+(p+1)\sum_{k=1}^{p}\binom{p}{k}\frac{H_{k}}{(e^{t}-1)^{k+1}},
\end{equation}
where $H_{n}$ is the harmonic numbers defined in \cite{Graham:1994:CMF:562056}:
\begin{equation*}
H_{n} = \sum_{j=1}^{n}\frac{1}{j}.
\end{equation*}
\end{theorem}

We provide a shorter proof of this theorem using elementary methods. 

\section{Proof}

For $p=0$ we have
\begin{equation*}
\sum_{n=0}^{\infty}B_{n, 0}\frac{t^{n}}{n!} = \sum_{n=0}^{\infty}B_{n}\frac{t^{n}}{n!} = \frac{t}{e^{t}-1}
\end{equation*}
and
\begin{equation*}
\frac{(0+1)(t-H_{0})e^{0}}{(e^{t}-1)^{0+1}}+(0+1)\sum_{k=1}^{0}\binom{0}{k}\frac{H_{k}}{(e^{t}-1)^{k+1}} = \frac{t}{e^{t}-1}.
\end{equation*}
Therefore Equation (\ref{eq:main}) holds for $p=0$. 

Now, assume that it holds for some $p$. Then 
\begin{equation*}
f_{p}(t) = \frac{(p+1)(t-H_{p})e^{pt}}{(e^{t}-1)^{p+1}}+(p+1)\sum_{k=1}^{p}\binom{p}{k}\frac{H_{k}}{(e^{t}-1)^{k+1}}.
\end{equation*}
It follows that

\begin{align*}
pf_p(t)&-f_p'(t)=p(p+1)\frac{(t-H_p)e^{pt}}{(e^{t}-1)^{p+1}}+p(p+1)\sum_{k=1}^p\binom{p}{k}\frac{H_k}{(e^{t}-1)^{k+1}}\\
&\quad -(p+1){\frac {{e^{pt}}}{  ( {e^{t}}-1  ) ^{p+1}}}-p(p+1){
	\frac { ( t-H_p  ) {e^{pt}}}{ ( {e^{t}}-1
		 ) ^{p+1}}}+(p+1)^2{\frac { ( t-H_p  ) e^{(p+1)t}}{ ( {e^{t}}-1  ) ^{p+2} }}
\\
&\quad +(p+1)\sum_{k=1}^p\binom{p}{k}(k+1)\frac{H_ke^t}{(e^{t}-1)^{k+2}}\\
&=p(p+1)\sum_{k=1}^p\binom{p}{k}\frac{H_k}{(e^{t}-1)^{k+1}}-(p+1){\frac {{e^{pt}}}{  ( {e^{t}}-1  ) ^{p+1}}}+(p+1)^2{\frac { ( t-H_p  ) {e^{(p+1)t}} }{ ( {e^{t}}-1  ) ^{p+2} }}
\\
&\quad +(p+1)\sum_{k=1}^p\binom{p}{k}(k+1)\frac{H_k}{(e^{t}-1)^{k+1}}
+(p+1)\sum_{k=1}^p\binom{p}{k}(k+1)\frac{H_k}{(e^{t}-1)^{k+2}}\\
&=-(p+1){\frac {{e^{pt}}}{  ( {e^{t}}-1  ) ^{p+1}}}+(p+1)^2{\frac { ( t-H_p  ) {e^{(p+1)t}} }{ ( {e^{t}}-1  ) ^{p+2} }}
+p(p+1)\sum_{k=1}^p\binom{p}{k}\frac{H_k}{(e^{t}-1)^{k+1}}\\
&\quad+(p+1)\sum_{k=1}^p\binom{p}{k}(k+1)\frac{H_k}{(e^{t}-1)^{k+1}}
+(p+1)\sum_{k=1}^{p+1}\binom{p}{k-1}k\frac{H_{k}-\frac1k}{(e^{t}-1)^{k+1}}\\
&=-(p+1){\frac {{e^{pt}}}{  ( {e^{t}}-1  ) ^{p+1}}}+(p+1)^2{\frac { ( t-H_p  ) e^{(p+1)t}}{ ( {e^{t}}-1  ) ^{p+2} }}
\\
&\quad +(p+1)\sum_{k=1}^p\frac{H_k}{(e^{t}-1)^{k+1}}\bigg[p\binom{p}{k}+\binom{p}{k}(k+1)+\binom{p}{k-1}k\bigg]\\
&\quad -(p+1)\sum_{k=1}^{p}\binom{p}{k-1}\frac{1}{(e^{t}-1)^{k+1}} + (p+1)^{2}\frac{H_{p}}{(e^{t}-1)^{p+2}}\\ 
&=-{\frac {(p+1){e^{pt}}}{  ( {e^{t}}-1  ) ^{p+1}}}+(p+1)^2{\frac { ( t-H_{p+1}  ) e^{(p+1)t}}{ ( {e^{t}}-1  ) ^{p+2} }}
+{\frac { (p+1) e^{(p+1)t}}{ ( {e^{t}}-1  ) ^{p+2} }}\\
&\quad +(p+1)^2\sum_{k=1}^{p+1}\frac{H_k}{(e^{t}-1)^{k+1}}\binom{p+1}{k}-(p+1)^2\frac{H_{p+1}}{(e^{t}-1)^{p+2}}\\
&\quad -\frac{(p+1)e^{pt}}{(e^t-1)^{p+2}} + (p+1)^{2}\frac{H_{p+1}}{(e^{t}-1)^{p+2}}\\ 
&=-{\frac {(p+1){e^{pt}}}{  ( {e^{t}}-1  ) ^{p+1}}}
+{\frac { (p+1) e^{(p+1)t}}{ ( {e^{t}}-1  ) ^{p+2} }}-\frac{(p+1)e^{pt}}{(e^t-1)^{p+2}}\\
&\quad +(p+1)^2\sum_{k=1}^{p+1}\frac{H_k}{(e^{t}-1)^{k+1}}\binom{p+1}{k}+(p+1)^2{\frac { ( t-H_{p+1}  ) e^{(p+1)t}}{ ( {e^{t}}-1  ) ^{p+2} }}\\
&=(p+1)^2\sum_{k=1}^{p+1}\frac{H_k}{(e^{t}-1)^{k+1}}\binom{p+1}{k}+(p+1)^2{\frac { ( t-H_{p+1}  ) e^{(p+1)t}}{ ( {e^{t}}-1  ) ^{p+2} }}\\
 &= \frac{(p+1)^{2}}{p+2}f_{p+1}(t).
\end{align*} 
Therefore Equation (\ref{eq:main}) holds for all $p\geq 0$ and the elementary proof of Theorem \ref{thm} is complete.

\bibliography{bernoulli_closed}

\end{document}